\documentclass[10pt]{amsart}

\usepackage{latexsym,amssymb,amsmath,graphics,xy}
\usepackage[dvips]{graphicx}

\def\ffi{\varphi}

\def\dst{\displaystyle}

\def\C{{\mathbb{C}}}

\def\R{{\mathbb{R}}}
\def\S{{\mathbb{S}}}
\def\T{{\mathbb{T}}}

\def\Z{{\mathbb{Z}}}

\newcommand{\norm}[1]{{\left\|{#1}\right\|}}
\newcommand{\ent}[1]{{\left[{#1}\right]}}
\newcommand{\abs}[1]{{\left|{#1}\right|}}
\newcommand{\scal}[1]{{\left\langle{#1}\right\rangle}}

\newenvironment{notation}[1][]{\vskip1pt\noindent\rm\textbf{Notation.}\ }{\rm\vskip1pt}

\newenvironment{definition}[1][]{\vskip3pt\noindent\sl\textbf{Definition.}\ }{\rm\vskip3pt}

\newtheorem{lemma}{Lemma}[section]
\newtheorem{proposition}[lemma]{Proposition}

\newtheorem{corollary}[lemma]{Corollary}

\begin{document}

\title{Nazarov's uncertainty principles in higher dimension}
\author{Philippe Jaming}
\address{Universit\'e d'Orl\'eans\\
Facult\'e des Sciences\\ 
MAPMO - F\'ed\'eration Denis Poisson\\ BP 6759\\ F 45067 Orl\'eans Cedex 2\\
France}
\email{Philippe.Jaming@univ-orleans.fr}

\begin{abstract}
In this paper we prove that there exists a constant $C$ such that, if $S,\Sigma$
are subsets of $\R^d$ of finite measure, then for every function $f\in L^2(\R^d)$,
$$
\int_{\R^d}|f(x)|^2\,\mbox{d}x\leq C e^{C\min\bigl(|S||\Sigma|,|S|^{1/d}w(\Sigma),w(S)|\Sigma|^{1/d}\bigr)}\left(
\int_{\R^d\setminus S}|f(x)|^2\,\mbox{d}x+\int_{\R^d\setminus\Sigma}|\widehat{f}(x)|^2\,\mbox{d}x\right)
$$
where $\widehat{f}$ is the Fourier transform of $f$ and $w(\Sigma)$ is the mean width of $\Sigma$. This extends 
to dimension $d\geq 1$ a result of Nazarov \cite{pp.Na} in dimension $d=1$.
\end{abstract}



\maketitle

\section{Introduction}

An uncertainty principle is a mathematical result that gives
limitations on the simultaneous localization of a function and its Fourier transform.
There are many statements of that nature, the most famous being due to Heisenberg-Pauli-Weil
when localization is measured in terms of smallness of dispersions 
and to Hardy when localization is measured in terms of fast decrease of the functions.
We refer the reader to the surveys \cite{pp.FS,pp.BD} and to the book \cite{pp.HJ}
for further references and results.

We will need a few notations before going on. In this paper $d$ will be a positive integer, all subsets of
$\R^d$ considered will be measurable and we will denote by $|S|$ the Lebesgue measure of $S$.
The Fourier transform is defined for $f\in L^1(\R^d)\cap L^2(\R^d)$ by
$$
\widehat{f}(\xi)=\int_{\R^d}f(x)e^{2i\pi\scal{x,\xi}}\,\mbox{d}x
$$
and extended to all of $L^2(\R^d)$ in the usual way.

In this paper, we are interested in another criterium of localization, namely smallness of support.
For instance, it is well known that if a function is compactly supported, then its Fourier transform
is an entire function and can therefore not be compactly supported. We may then ask what happens
if a function $f$ and its Fourier transform $\widehat{f}$ are only small outside a compact set? This leads naturally
to the following definition:

\medskip

\noindent{\bf Definition.}\\
{\sl Let $S,\Sigma$ be two Borel subsets of $\R^d$. Then we will say that

--- $(S,\Sigma)$ is an annihilating pair (a-pair in short) if the only function $f$ that is supported in $S$
and such that its Fourier transform $\widehat{f}$ is supported in $\Sigma$ is $f=0$;

--- $(S,\Sigma)$ is a strong annihilating pair (strong a-pair in short) if there exists a constant $C=C(S,\Sigma)$
such that for every $f\in L^2(\R^d)$,}
$$
\int_{\R^d}|f(x)|^2\,\mbox{d}x\leq C\left(
\int_{\R^d\setminus S}|f(x)|^2\,\mbox{d}x+\int_{\R^d\setminus\Sigma}|\widehat{f}(x)|^2\,\mbox{d}x\right).
$$

\medskip

This notion has been extensively studied in the case $S$ is a compact set by Logvinenko and Sereda \cite{pp.LS},
Paneah \cite{pp.Pa1,pp.Pa2}, Havin and J\"oricke \cite{pp.HJ2} and Kovrijkine \cite{pp.Ko}, see
also \cite{pp.HJ}. In this case the class of all
$\Sigma$'s for which $(S,\Sigma)$ is a strong a-pair is characterized. Moreover, if $S$ is convex,
there are fairly good estimates of the constant $C(S,\Sigma)$ in terms of the geometry of $S$
and $\Sigma$. 

For sets $S,\Sigma$ that are sublevel sets of quadratic forms, the problem has been studied by
Shubin, Vakilian, Wolff \cite{pp.SVW} and by Demange \cite{pp.De1,pp.De2}.

Here we will focus on the case of $S,\Sigma$ being of finite Lebesgue measure. This was first studied
by Benedicks \cite{pp.Be} who proved that in this case $(S,\Sigma)$ is an a-pair,
and a little abstract nonsense allows to prove that in this case $(S,\Sigma)$ is also a strong a-pair,
see \cite{pp.BD}. This last fact was proved with a different method by Amrein and Berthier \cite{pp.AB}.
Unfortunatly both proofs do not give any estimate on the constant $C(S,\Sigma)$. By using
a randomization of Benedicks proof and an extension of a lemma of Turan, Nazarov \cite{pp.Na}
showed that in dimension $1$, the constant is of the form $C(S,\Sigma)=Ce^{C|S||\Sigma|}$.
It was thought for some time that Nazarov's method would extend to higher dimension
to give a constant of the same form. This is far from the expected optimal which is thought to be obtained
by taking $S,\Sigma$ balls of radius $R$ and $f$ a Gaussian function, which gives
$C(S,\Sigma)=Ce^{CR^2}=e^{C(|S||\Sigma|)^{1/d}}$.

The aim of this paper is to push Nazarov's technique as far as possible and thus improve
the $Ce^{C|S||\Sigma|}$ constant when the geometry of $\Sigma$ is suitable. Using the recent extension
of Nazarov's Turan lemma to higher dimension by Fontes-Merz \cite{pp.FM}, we will prove the following result:

\medskip

\noindent{\bf Theorem.}\\
{\sl There exists a constant $C$ such that, for every sets $S,\Sigma\subset\R^d$ of finite
Lebesgue measure and for every $f\in L^2(\R^d)$,
$$
\int_{\R^d}|f(x)|^2\,\mbox{d}x\leq C e^{C\min\bigl(|S||\Sigma|,|S|^{1/d}w(\Sigma),w(S)|\Sigma|^{1/d}\bigr)}\left(
\int_{\R^d\setminus S}|f(x)|^2\,\mbox{d}x+\int_{\R^d\setminus\Sigma}|\widehat{f}(x)|^2\,\mbox{d}x\right)
$$
where $w(\Sigma)$ is the mean width of $\Sigma$.}

\medskip

In particular, if $S$ or $\Sigma$ has a geometry that is close to a ball, this is in accordance with what is supposed to be the optimal result.

\medskip

The remaining of this paper is devoted to the proof of this theorem. In order to do so, we first extend to higher dimension the random periodization technique. Then we recall the Turan type estimates we will need. The last section is then devoted to the proof
of the theorem.

\section{Random Periodization}

\subsection{Preliminaries}\ \\
For any integer $d$, let $SO(d)$ denote the group of rotations on $\R^d$.
Denote by $\mbox{d}\nu_d$ the normalized Haar measure on $SO(d)$.
Then there exists a constant $C=C(d)$ such that, for every $u\in\S^{d-1}$, the unit sphere $\S^{d-1}$ of $\R^d$,
and every function $f\in L^1(\R^d)$
$$
\int_{SO(d)}\int_0^{+\infty} f\bigl(v\,\rho(u)\bigr)v^{d-1}\,\mbox{d}v\,\mbox{d}\nu_d(\rho)=
C\int_{\R^d}f(x)\,\mbox{d}x.
$$

\subsection{The higher dimensional Lattice Averaging Lemma}\ \\
The following lemma was proved by Nazarov in dimension $d=1$.

\begin{lemma}[Lattice Averaging Lemma]
\label{inc.lem:LAL}\ \\
Let $d\geq 1$ be an integer, then 
for every $\ffi\in L^1(\R^d)$, $\ffi\geq0$, the following estimates hold
$$
\int_{SO(d)}\int_1^2\sum_{k\in\Z^d\setminus\{0\}} \ffi\bigl(v\,\rho(k)\bigr)\,\mathrm{d}v\,\mathrm{d}\nu_d(\rho)
\simeq \int_{\|x\|\geq 1}\ffi(x)\,\mathrm{d}x
$$
and
\begin{equation}
\label{naz.eq.int}
\int_{SO(d)}\int_1^2\sum_{k\in\Z^d\setminus\{0\}} \ffi\left(\frac{\rho(k)}{v}\right)\,\mathrm{d}v\,\mathrm{d}\nu_d(\rho)
\simeq \int_{\|x\|\geq 1/2}\ffi(x)\,\mathrm{d}x.
\end{equation}
\end{lemma}

Here, as usual, by $A\simeq B$ we mean that there exists a constant $C$ depending only on $d$
such that $\frac{1}{C}B\leq A\leq CB$.

\begin{proof} With \eqref{naz.eq.int}, we get
\begin{eqnarray*}
\int_{SO(d)}\int_1^2\sum_{k\in\Z^d\setminus\{0\}}\ffi\bigl(v\,\rho(k)\bigr)\mbox{d}v\,\mbox{d}\nu_d(\rho)
&\simeq&\sum_{k\in\Z^d\setminus\{0\}}
\int_{SO(d)}\int_1^2 \ffi\bigl(v\norm{k}\,\rho(k/\norm{k})\bigr)\,v^{d-1}\mbox{d}v\,\mbox{d}\nu_d(\rho)\\
&=&C\sum_{k\in\Z^d\setminus\{0\}}\int_{1\leq\norm{x}\leq 2}\ffi(\norm{k}x)\,\mbox{d}x\\
&=&C\sum_{k\in\Z^d\setminus\{0\}}\frac{1}{\|k\|^d}
\int_{\|k\|\leq \|x\|\leq 2\|k\|}\ffi(x)\,\mbox{d}x\\
&=&C\int_{\|x\|\geq 1} \ffi(x)\,\sum_{\|k\|\leq\|x\|\leq 2\|k\|}\frac{1}{\|k\|^d}\,\mbox{d}x\\
&\simeq&\int_{\|x\|\geq 1} \ffi(x)\,\mbox{d}x
\end{eqnarray*}
since, for $\norm{x}\geq 1$,
$$
\sum_{\|k\|\leq\|x\|\leq 2\|k\|}\frac{1}{\|k\|^d}\simeq
\frac{\abs{\left\{u\in\R^d\,:\ \|x\|/2\leq u\leq \|x\|\right\}}}{\|x\|^d}
=|B(0,1)\setminus B(0,1/2)|.
$$

For the second statement, one first changes $v$ into $1/v$ and the remaining of the proof is similar.
\end{proof}

\begin{definition}\ \\
For a function $f\in L^2(\R)$, $\rho\in SO(d)$ and $v>0$,
we define the periodization
$\Gamma_{\rho,v}(t)=\Gamma_{\rho,v}(f)(t)$ of the function $f$ by
$$
\Gamma_{\rho,v}(t)= \frac{1}{\sqrt{v}}\sum_{k\in\Z^{d}}f\left(\frac{\rho(k+t)}{v}\right).
$$
\end{definition}

The series in the definition of $\Gamma_{\rho,v}$ converges in $L^2(\T^d)$ and represents a periodic function.
An easy computation shows that the Fourier coefficients of $\Gamma_{\rho,v}$ are
$\widehat{\Gamma_{\rho,v}}(m)=\sqrt{v}\widehat{f}\bigl(v\,{}^t\rho(m)\bigr)$ for $m\in\Z^d$.

\begin{notation}\ \\
In the sequel, $v$ will be considered as a random variable equidistributed on the interval $(1,2)$
and $\rho$ as a random variable equidistributed on $SO(d)$. 
The expectation with respect to these random variables will be denoted by $\mathbb{E}_{\rho,v}$
\end{notation}

\subsection{Properties of random periodizations}\ \\
From the Lattice Averaging Lemma we shall derive the following simple but useful properties
of the random periodization.

\begin{proposition}
\label{inc.prop:conslala}\ \\
Let $d\geq 1$ be an integer and $C=C(d)$ be the constant defined in Lemma \ref{inc.lem:LAL}. 
Let $S\subset\R^d$ be a set of finite measure and let $f\in L^2(\R^d)$ be supported in $S$. Then
\begin{enumerate}
\renewcommand{\theenumi}{\roman{enumi}}
\item for all $v\in(1,2)$, $|\{\,t\in(0,1)\,:\ \Gamma_{\rho,v}(t)\ne 0\,\}\leq 2^d|S|$;

\item $\mathbb{E}_{\rho,v}\bigl(\|\Gamma_{\rho,v}\|_{L^2(0,1)}^2\bigr)
\leq 2 |\widehat{f}(0)|^2+2C\|f\|_{L^2(\R^d)}^2
\leq 2(|S|+C)\|f\|_{L^2(\R^d)}^2$.
\end{enumerate}
\end{proposition}

\begin{proof} $i)$ The set of all points $t\in[0,1]^d$ for which the summand
$f\left(\frac{\rho(k+t)}{v}\right)$ in the series defining $\Gamma_{\rho,v}$ does not vanish equals
$v\,{}^t\rho(S)\cap\bigl([0,1]^d+k\bigr)$. Therefore,
$$
|\{t\in[0,1]^d\,:\ \Gamma_{\rho,v}(t)\neq 0\}|\leq \sum_{k\in\Z^d}|v\,{}^t\rho(S)\cap\bigl([0,1]^d+k\bigr)|
=|v\,{}^t\rho(S)|\le 2^d|S|.
$$
$ii)$ Parseval's Identity gives
$$
\mathbb{E}_{\rho,v}\bigl(\|\Gamma_{\rho,v}\|_{L^2(\T^d)}^2\bigr)
=\mathbb{E}_{\rho,v}\left(\sum_{k\in\Z^d}|\widehat{\Gamma_{\rho,v}}(k)|^2\right)
=\mathbb{E}_{\rho,v}\bigl(|\widehat{\Gamma_{\rho,v}}(0)|^2\bigr)
+\mathbb{E}_{\rho,v}\left(\sum_{k\in\Z^d\setminus\{0\}}|\widehat{\Gamma_{\rho,v}}(k)|^2\right).
$$
But
$|\widehat{\Gamma_{\rho,v}}(0)|^2= v|\widehat{f}(0)|^2 \le 2|\widehat{f}(0)|^2$, and,
with the Lattice Averaging Lemma,
\begin{eqnarray*}
\mathbb{E}_{\rho,v}\left(\sum_{m\in\Z^d\setminus\{0\}}|\Gamma_{\rho,v}(m)|^2\right)
&=&\int_{SO(d)}\int_1^2\left(\sum_{m\in\Z^d\setminus\{0\}} v|\widehat{f}\bigl(v\,\rho(m)\bigr)|^2\right)\,\mbox{d}v\,\mbox{d}\nu_{d}(\rho)\\
&\leq&2\int_{SO(d)}\int_1^2\left(\sum_{m\in\Z^d\setminus\{0\}}|\widehat{f}\bigl(v\,\rho\zeta(m)\bigr)|^2 \right)\,\mbox{d}v\,\mbox{d}\nu_{d}(\rho)\\
&\leq& 2C\int_{\R^d}|\widehat{f}\bigl(\rho(\xi)\bigr)|^2\mbox{d}\xi=2C\|f\|_{L^2(\R^d)}^2.
\end{eqnarray*}
It remains to notice that
$$
|\widehat{f}(0)|^2=\abs{\int_S f(x)\mbox{d}x}^2\leq |S|\int_S|f(x)|^2\mbox{d}x
=|S|\|f\|_{L^2(\R)}^2.
$$
\end{proof}

\begin{definition}\ \\
Let $\Sigma\subset\R$ be a measurable set with, $0\in\Sigma$. We consider the lattice
$\Lambda=\Lambda(\rho,v):=\{v\,{}^t\rho(j)\,:\ j\in\Z^d\}$
and denote $\mathcal{M}_{\rho,v}=\{k\in\Z^d\,:\ v\,{}^t\rho(k)\in\Sigma\}=\Lambda\cap\Sigma$.
\end{definition}

\begin{proposition}
\label{inc.prop:conslalb}\ \\
With the previous notations
\begin{enumerate}
\renewcommand{\theenumi}{\roman{enumi}}
\item $\dst\mathbb{E}_{\rho,v}\bigl(\mathrm{card}\,\mathcal{M}_{\rho,v}-1\bigr)\leq C|\Sigma|$, in particular $\mathcal{M}_{\rho,v}$ is almost surely finite;

\item $\dst\mathbb{E}_{\rho,v}\left(\sum_{m\in\Z^d\setminus\mathcal{M}_{\rho,v}}|\widehat{\Gamma_{\rho,v}}(m)|^2\right)
 \leq 2C\int_{\R^d\setminus\Sigma}|\widehat{f}(\xi)|^2\,\mathrm{d}\xi$.
\end{enumerate}
\end{proposition}

\begin{proof}
$i)$ Since $\mathrm{card}\,\mathcal{M}_{\rho,v}=1+\sum_{m\in\Z^d\setminus\{0\}}\chi_\Sigma\bigl(v\,{}^t\rho(m)\bigr)$,
we have
\begin{eqnarray*}
\mathbb{E}_{\rho,v}\bigl(\mathrm{card}\,\mathcal{M}_{\rho,v}-1\bigr)
&=&\int_{SO(d)}\int_1^2 \sum_{k\in\Z^d\setminus\{0\}}\chi_\Sigma\bigl(v\,{}^t\rho(k)\bigr)\,\mbox{d}v\,\mbox{d}\nu_{d}(\rho)\\
&\leq& C\int_{\R^d}\chi_\Sigma(x)\,\mbox{d}x=C|\Sigma|.
\end{eqnarray*}
$ii)$ From the expression of $\widehat{\Gamma_{\rho,v}}$ we get that 
$\dst \mathbb{E}_{\rho,v}\left(\sum_{m\in\Z^d\setminus\mathcal{M}_{\rho,v}}|\widehat{\Gamma_{\rho,v}}(k)|^2\right)$
is
\begin{eqnarray*}
&=&\int_{SO(d)}\int_1^2\left(\sum_{m\in\Z^d\setminus\{0\}} v\bigl|\widehat{f}\bigl(v\,{}^t\rho(m)\bigr)\bigr|^2
\chi_{\R^d\setminus\Sigma}\bigl(v\,{}^t\rho(m)\bigr)\right)\,\mbox{d}v\,\mbox{d}\nu_d(\rho)\\
&\leq&2\int_{SO(d)}\int_1^2\left(\sum_{m\in\Z^d\setminus\{0\}}
\bigl|\widehat{f}\bigl(v\,{}^t\rho(mk)\bigr)\bigr|^2\chi_{\R^d\setminus\Sigma}\bigl(v\,{}^t\rho(m)\bigr)\right)
\,\mbox{d}v\,\mbox{d}\nu_d(\rho) \\
&\leq& 2C\int_{\R^d}|\widehat{f}(\xi)|^2\chi_{\R\setminus\Sigma}(\xi)\,\mbox{d}\xi
=2C\int_{\R^d\setminus\Sigma}|\widehat{f}(\xi)|^2\,\mbox{d}\xi.
\end{eqnarray*}
by Lemma \ref{inc.lem:LAL}.
\end{proof}

\section{A Turan Lemma}

\subsection{Nazarov and Fontes-Merz' Turan Lemmas}\ \\
For sake of completeness, we will recall here the Turan type estimates of
trigonometric polynomials we will need.

\medskip

\noindent{\bf Theorem (Nazarov's Turan Lemma \cite{pp.Na})}\\
{\sl Let $\dst P(t)=\sum_{k=1}^m c_ke^{2i\pi r_k t}$ with $c_k\in\C\setminus\{0\}$, $r_1<\dots<r_m\in\Z$, be a
trigonometric polynomial of \emph{order} $\mathrm{ord}\,P=m$
and let $E$ be a measurable subset of $\T$. Then}
\begin{equation}
\label{eq:turannaz}
\sup_{z\in \T}|P(z)|\leq\left(\frac{14}{|E|}\right)^{m-1}\sup_{z\in E}|P(z)|.
\end{equation}

\medskip

The original theorem of Turan deals with sets $E$ that are arcs.
The extension to higher dimension has been obtained in \cite{pp.FM} using a clever induction on
the dimension. 

\medskip

\noindent{\bf Corollary (Fontes-Merz's Turan Lemma \cite{pp.FM})}\\
{\sl 
Let $d\geq 1$ be an integer and let
$$
p(z_1,\ldots,z_d)=\sum_{k_1=0}^{m_1}\cdots\sum_{k_d=0}^{m_d}c_{k_1,\ldots,k_d}z_1^{r_{1,k_1}}\cdots z_d^{r_{r_d,k_d}}
$$
with $r_{i,k_i}\in\Z$ be a polynomial in $d$ variables. Then, for every measurable set $E\subset\T^d$,}
$$
\sup_{z\in\T^d}|p(z)|\leq\left(\frac{14d}{|E|}\right)^{m_1+\cdots+m_d}
\sup_{z\in E}|p(z)|.
$$

\medskip

The quantity $m_1+\cdots+m_d$ is called the order of $p$ (with the usual convention that we take
the most compact possible representation of $p$) and is denoted by $\mbox{ord}\,p$. In general
$$
m_1+\cdots+m_d\leq d\max m_i\leq d\,\mbox{Card}\,\mbox{Spec}\,p
$$
while
$$
\mbox{Card}\,\mbox{Spec}\,p\leq (m_1+1)\cdots(m_d+1).
$$

\subsection{An estimate of the average order}\ \\
The notion of order of a polynomial suggests the following definition of the order
of a subset of $\Z^d$.

\noindent{\bf Definition.}\\
{\sl Let $M\subset\Z^d$ be a finite set, we will say that $M$ is of order $k$ 
and write $\mbox{ord}\,M=k$ if there exists integers
$m_1,\ldots,m_d$ with $m_1+\cdots+m_d=k$ such that the projection of
$M$ on the $i$-th coordinate axis has $m_i$ elements.

Finally, if $\Lambda=A\Z^d$ is a lattice and $M\subset \Lambda$ is finite, we will call $\mbox{ord}\,M=
\mbox{ord}\,A^{-1}M$.}

\medskip

Note that
$$
m_1=\sum_{k\in\Z}\sup_{k'\in\Z^{d-1}}\chi_{M}(k,k')
$$
with similar expressions for the other $m_i$'s.

In order to estimate the order of the set $\mathcal{M}_{\rho,v}$
introduced before Proposition \ref{inc.prop:conslalb}, the easiest is to bound
the order by the cardinal of the set, which amounts to bounding the supremum by the sum
over $k'\in\Z^{d-1}$ in the above expression. One then gets 
$\dst\mathbb{E}_{\rho,v}\bigl(\mathrm{ord}\,\mathcal{M}_{\rho,v}-d\bigr)\leq C|\Sigma|$.
This shows in particular that it is enough to estimate this quantity
when $\Sigma$ is a relatively compact open set.

The proof of the uncertainty principle in the next section will then give a constant
$Ce^{C|S||\Sigma|}$ in Nazarov's result. We will slightly improve this.
in order to do so, let us introduce the following quantities:

--- the average width: for a relatively compact open set $\Sigma$ 	and for $\rho\in SO(d)$,
let $P_\rho(\Sigma)$ be the projection of $\Sigma$ on the span of $\rho(1,0,\ldots,0)$.
We define
$$
w(\Sigma)=\int_{SO(d)}|P_\rho(\Sigma)|\,\mbox{d}\nu_d(\rho)
$$
the average width of $\Sigma$. If $\Sigma$ is a ball, this is just its diameter.

--- let us also introduce the measure $\mu$ on
$\R^d$ defined by
$$
\mu(\Sigma)
=\inf\left\{\sum_{i\in I}\min(r_i,r_i^d)\,:
\{B(x_i,r_i)\}_{i\in I}\mbox{ is a cover of }\Sigma\right\}.
$$
Note that $\mu(\Sigma)\leq C|\Sigma|$ since the $d$-dimensional Hausdorff measure is the Lebesgue measure.

We will now prove the following:

\begin{proposition}
\label{pp.prop:ord}\ \\
Let $\Sigma$ be a relatively compact open set with $0\in\Sigma$. We consider a random lattice
$\Lambda=\Lambda(\rho,v):=\{v\,{}^t\rho(j)\,:\ j\in\Z^d\}$
and denote $\mathcal{M}_{\rho,v}=\{k\in\Z^d\,:\ v\,{}^t\rho(k)\in\Sigma\}=\Lambda\cap\Sigma$.
Then $\dst\mathbb{E}_{\rho,v}\bigl(\mathrm{ord}\,\mathcal{M}_{\rho,v}-d\bigr)\leq C\min\bigl(\mu(\Sigma),w(\Sigma)\bigr)$.
\end{proposition}

\begin{proof} Let
$$
m_{\rho,v}(\Sigma)=\sum_{k\in\Z\setminus\{0\}}\sup_{k'\in\Z^{d-1}}\chi_{\Sigma}\bigl(v{}^t\rho(k,k')\bigr).
$$
It is enough to prove that
\begin{equation}
\label{pp.eq:enough}
\mathbb{E}_{\rho,v}\bigl(m_{\rho,v}(\Sigma)\bigr)\leq C\min\bigl(\mu(\Sigma),w(\Sigma)\bigr).
\end{equation}
As pointed out above, $\mathbb{E}_{\rho,v}\bigl(m_{\rho,v}(\Sigma)\bigr)\leq C|\Sigma|$. In particular,
if $\Sigma$ is a ball of radius $r$, $\mathbb{E}_{\rho,v}\bigl(m_{\rho,v}(\Sigma)\bigr)\leq Cr^d$

On the other hand
$$
m_{\rho,v}:=m_{\rho,v}\bigl(\Sigma\bigr)\leq\sum_{k\in\Z\setminus\{0\}}\sup_{y\in\R^{d-1}}\chi_{\Sigma}\bigl({}^t\rho(vk,y)\bigr)
$$
and the one-dimensional lattice averaging lemma then gives
\begin{eqnarray*}
\mathbb{E}_{\rho,v}(m_{\rho,v})&\leq& \int_{SO(d)}\int_{|x|\geq 1}\sup_{y\in\R^{d-1}}\chi_{B\bigl(\rho(a),r\bigr)}(x,y)
\,\mbox{d}x\,\mbox{d}\nu_d(\rho)\\
&\leq&C\int_{SO(d)}\int_{|x|\geq 1}\chi_{P_{{}^t\rho}}(x)\,\mbox{d}x\,\mbox{d}\nu_d(\rho)\\
&\leq&C w(\Sigma).
\end{eqnarray*}

In particular, if $\Sigma$ is a ball of radius $r$, then $\mathbb{E}_{\rho,v}\bigl(m_{\rho,v}(\Sigma)\bigr)\leq Cr$.
To conclude, it is enough to
note that $m_{\rho,v}(\Sigma\cup\Sigma')\leq m_{\rho,v}(\Sigma)+m_{\rho,v}(\Sigma')$ and that
if $\Sigma\subset\Sigma'$ then $m_{\rho,v}(\Sigma)\leq m_{\rho,v}(\Sigma')$. Covering
$\Sigma$ with balls then gives the desired result.
\end{proof}

The result above is essentially sharp as the following example shows.
For simplicity, we will give the example in dimension $d=2$. Let $N\geq1$ be an integer and
let $R\gg 1$ be two real numbers. Let
$$
\Sigma_N=\bigcup_{j=0}^{N-1}B\left(R(\cos\frac{2\pi}{N}j,\sin\frac{2\pi}{N}j),\frac{1}{2}\right).
$$
That is, $\Sigma_N$ is the union of $N$ discs regularily placed on a big circle,
see the figure below. 
\begin{center}
\begin{tabular}{c}
\begin{xy}
<18mm,-22mm>*{} ; <18mm,22mm>*{} **\dir{.},
<0mm,0mm>*{} ; <20mm,25mm>*{} **\dir{-},
<8mm,10mm>*{} ;<-7mm,22mm>*{} **\dir{.},
<8mm,10mm>*{} ;<23mm,-2mm>*{} **\dir{.},
<20mm,0mm>*\cir<2mm>{},
<14mm,14mm>*\cir<2mm>{},
<0mm,20mm>*\cir<2mm>{},
<14mm,-14mm>*\cir<2mm>{},
<0mm,-20mm>*\cir<2mm>{},
<-20mm,0mm>*\cir<2mm>{},
<-14mm,14mm>*\cir<2mm>{},
<0mm,-20mm>*\cir<2mm>{},
<-14mm,-14mm>*\cir<2mm>{},
\end{xy}\\\ \\
\begin{minipage}{12.5cm}
\begin{center}
The set $\Sigma_N$~:\\
\end{center}
Note that each line orthogonal to a line through the origin meets
at most two circles. Moreover, these circles have radius $\leq 1/2$
thus, for $k$ fixed, at most two segments $\{{}^t\rho(vk,vk'),v\in(1,2)\}$
can intersect $\Sigma_N$. Therefore, the $\sup_{k'}$ in the formula defining $m_1$
can be bounded below by $\frac{1}{2}\sum_{k'}$.
\end{minipage}
\end{tabular}
\end{center}

\medskip

Then, for each $k$, $\#\{k'\in\Z\,:\ {}^t\rho(vk,vk')\cap\Sigma_N\not=0\}\leq 2$
thus
$$
m_{\rho,v}\bigl(\Sigma_N\bigr)\geq \frac{1}{2}\sum_{k\in\Z\setminus\{0\}}\sum_{k'\in\Z}\chi_{\Sigma_N}\bigl({}^t\rho(vk,vk')\bigr)
\simeq |\Sigma_N|\simeq N\simeq w(\Sigma_N)
$$
with the Latice Averaging Lemma.

\section{Conclusion}
The remaining of the proof follows the path of Nazarov's original argument. We
include it here for sake of completeness.

Let us write $\nu(\Sigma)=\min\bigl(w(\Sigma),\mu(\Sigma)\bigr)$.
First, it is enough to prove that there exists a constant $C=C(d)$ such that
$$
\int_{\Sigma}|\widehat{f}(\xi)|^2\,\mbox{d}\xi
\leq Ce^{C\nu(|S|^{1/d}\Sigma)}\int_{\R^d\setminus\Sigma}|\widehat{f}(\xi)|^2\,\mbox{d}\xi
$$
for every $f\in L^2(S)$.
Moreover, using a scaling argument, it is enough to show that,
if $|S|=2^{-d+1}$, then for every set $\Sigma$ and every $f\in L^2(S)$,
$$
\int_{\Sigma}|\widehat{f}(\xi)|^2\,\mbox{d}\xi
\leq Ce^{C\mu(\Sigma)}\int_{\R^d\setminus\Sigma}|\widehat{f}(\xi)|^2\,\mbox{d}\xi.
$$

Set $\Gamma_{\rho,v}(t)=\Gamma_{\rho,v}(f)(t)$ the random periodization of $f$.
Then, setting $E_{\rho,v}=\{t\in(0,1)\,:\ \Gamma_{\rho,v}(t)=0\}$,
we have by Proposition \ref{inc.prop:conslala} $i)$ that $|E_{\rho,v}|\geq 1-2^d|S|=\frac{1}{2}$.

Next, set $\mathcal{M}_{\rho,v}:=\{m\in\Z^d\,:\ v\,{}^t\rho(m) \in\Sigma\cup\{0\}\}$ and decompose
$\Gamma_{\rho,v}=P_{\rho,v}+R_{\rho,v}$ where
$$
P_{\rho,v}(t)=\sum_{m\in\mathcal{M}_{\rho,v}}\widehat{\Gamma_{\rho,v}}(m)e^{2i\pi mt}
$$
while
$$
R_{\rho,v}(t)=\sum_{m\in\Z^d\setminus\mathcal{M}_{\rho,v}}\widehat{\Gamma_{\rho,v}}(m)e^{2i\pi mt}.
$$

By Proposition \ref{inc.prop:conslalb} $ii)$,
$$
\mathbb{E}_{\rho,v}\bigl(\norm{R_{\rho,v}}_{L^2(0,1)}^2\bigr)=
\mathbb{E}_{\rho,v}\left(\sum_{m\in\Z^d\setminus\mathcal{M}_{\rho,v}}|\widehat{\Gamma_{\rho,v}}(m)|^2\right)
\leq 2C\int_{\R^d\setminus\Sigma}|\widehat{f}(\xi)|^2\,\mbox{d}\xi,
$$
hence
$$
\mathbb{E}_{\rho,v}\left(\norm{R_{\rho,v}}_{L^2(0,1)}^2>4C\int_{\R^d\setminus\Sigma}|\widehat{f}(\xi)|^2\,\mbox{d}\xi\right)
<\frac{1}{2}.
$$
On the other hand, by  Proposition \ref{pp.prop:ord},
$$
\mathbb{E}_{\rho,v}\bigl(\mathrm{ord}\,P_{\rho,v}\bigr)\leq C\mu(\Sigma)+d
$$
and therefore
$$
\mathbb{E}_{\rho,v}\bigl(\mathrm{ord}\,P_{\rho,v}>2\bigl(C\mu(\Sigma)+d)\bigr)<\frac{1}{2}.
$$
We thus get that the two events
\begin{enumerate}
\item $\dst\norm{R_{\rho,v}}_{L^2(0,1)}^2\leq 4C\int_{\R^d\setminus\Sigma}|\widehat{f}(\xi)|^2\,\mbox{d}\xi$,

\item $\mathrm{ord}\,P_{\rho,v}\leq 2\bigl(C\mu(\Sigma)+d)$
\end{enumerate}
happen simultaneously with non-zero probability, while the two events
\begin{enumerate}
\setcounter{enumi}{2}
\item $\dst|E_{\rho,v}|\geq\frac{1}{2}$,

\item $|\widehat{f}(0)|^2\leq |\widehat{P_{\rho,v}}(0)|^2$
\end{enumerate}
are certain. We will now take $v\in(1,2)$, $\rho\in\S^{d-1}$ such that all four events hold simultaneously.

Further, by definition $\Gamma_{\rho,v}=0$ on $E_{\rho,v}$, that is $P_{\rho,v}$ and $-R_{\rho,v}$ coincide on $E_{\rho,v}$.
It follows that
$$
\int_{E_{\rho,v}}|P_{\rho,v}(x)|^2\,\mbox{d}x
=\int_{E_{\rho,v}}|R_{\rho,v}(x)|^2\,\mbox{d}x
.$$
Hence
$$
\abs{\left\{x\in E_{\rho,v}\,: |P_{\rho,v}(x)|^2\geq 
16C\int_{\R^d\setminus\Sigma}|\widehat{f}(\xi)|^2\,\mbox{d}\xi\right\}}
\leq\frac{1}{4}
$$
and, as $|E_{\rho,v}|\geq\frac{1}{2}$, we get that $|\tilde E_{\rho,v}|\geq\frac{1}{4}$ where
$$
\tilde E_{\rho,v}=\left\{x\in E_{\rho,v}\,: |P_{\rho,v}(x)|\leq 
4\left(C\int_{\R^d\setminus\Sigma}|\widehat{f}(\xi)|^2
\,\mbox{d}\xi\right)^{\frac{1}{2}}\right\}.
$$

We can now apply Turan's Lemma and get
\begin{eqnarray*}
|\widehat{f}(0)|^2&\leq&|\widehat{P_{\rho,v}}(0)|^2\leq
\left(\sum_{k\in\Z^d}|\widehat{P_{\rho,v}}(k)|\right)^{2}\leq\bigl(\sup_{x\in\T^d}|P_{\rho,v}(x)|\bigr)^2\\
&\leq&\ent{\left(\frac{14d}{|\tilde E_{\rho,v}|}\right)^{\mathrm{ord}P_{\rho,v}-1}\sup_{x\in\tilde E_{\rho,v}}|P_{\rho,v}(x)|}^2\\
&\leq&\ent{\left(\frac{14d}{1/4}\right)^{\mathrm{ord}P_{\rho,v}-1} 4
\left(C\int_{\R^d\setminus\Sigma}|\widehat{f}(\xi)|^2\,\mbox{d}\xi\right)^{1/2}}^2\\
&\leq& Ce^{C\nu(\Sigma)}\int_{\R^d\setminus\Sigma}|\widehat{f}(\xi)|^2\,\mbox{d}\xi.
\end{eqnarray*}

If we now apply this to $f_y(x)=f(x)e^{-2i\pi xy}$ instead of $f$ and to the set $\Sigma_y=\Sigma-y$ instead of $\Sigma$,
we obtain that
$$
|\widehat{f}(y)|^2\leq C e^{C\nu(\Sigma)}\int_{\R^d\setminus\Sigma}|\widehat{f}(\xi)|^2\,\mbox{d}\xi.
$$
and integrating this over $\Sigma$ gives
$$
\int_\Sigma|\widehat{f}(y)|^2\,\mbox{d}y\leq C|\Sigma| e^{C\mu(\Sigma)}\int_{\R^d\setminus\Sigma}|\widehat{f}(\xi)|^2\,\mbox{d}\xi
$$
as claimed. \hfill$\Box$

The values of the constants may be tracked and linked to those of the Random Averaging Lemma, but
we do not expect these constants to be any near to optimal (as they are already not optimal in dimension 1)
so we will not pursue this.

Note also that, with mutadis mutandis the same proof as in \cite{pp.Na} we obtain the following corollary:

\begin{corollary}\ \\
Let $S,\Sigma$ be two measurable subsets of $\R^d$ and let $C$ be the constant of the main theorem.
Then, for every $p\in(0,2)$ and every $f\in L^p(\R^d)$ with spectrum in $\Sigma$,
$$
\norm{f}_{L^p}^p\leq Ce^{Cp|S||\Sigma|}\int_{\R^d\setminus S}|f(x)|^p\,\mbox{d}x.
$$
\end{corollary}

\section*{Aknowledgements}
The author wishes to thank B. Demange, L. Grafakos for valuable conversations concerning the content of this
paper and A. Bonami for a careful reading of the manuscript.

This work was supported by the Hungarian-French Scientific and
Technological Governmental Cooperation, no. F-10/04

\end{document}